\documentclass{amsart}

\usepackage{fancyhdr}
\usepackage{lastpage}
\usepackage{stmaryrd,yhmath}

\newcommand{\bdis}{\begin{displaymath}}
\newcommand{\edis}{\end{displaymath}}
\newcommand{\be}{\begin{equation}}
\newcommand{\ee}{\end{equation}}
\newcommand{\mbb}{\mathbb}
\newcommand{\mcal}{\mathcal}

\newcommand{\vp}{\varphi}

\newcommand{\zf}{\zeta\left(\frac{1}{2}+it\right)}

\theoremstyle{definition}

\newtheorem{cor}[]{Corollary}

\theoremstyle{remark}
\newtheorem{remark}[]{Remark}

\newtheorem*{mydef11}{{\bf Theorem 1}}

\newtheorem*{mydef12}{{\bf Theorem 2}}

\newtheorem*{mydef3}{{\bf Problem}}

\newtheorem*{mydef5}{{\bf Lemma}}

\numberwithin{equation}{section}

\begin{document}

\title{Jacob's ladders and the $\tilde{Z}^2$-transformation of the orthogonal system of trigonometric functions}

\author{Jan Moser}

\address{Department of Mathematical Analysis and Numerical Mathematics, Comenius University, Mlynska Dolina M105, 842 48 Bratislava, SLOVAKIA}

\email{jan.mozer@fmph.uniba.sk}


\begin{abstract}
It is shown in this paper that there is a continuum set of orthogonal systems relative to the weight function $\tilde{Z}^2(t)$. The corresponding
integrals cannot be obtained in known theories of Balasubramanian, Heath-Brown and Ivic.
\end{abstract}

\maketitle

\section{The first result}

\subsection{}

In this paper we obtain some new properties of the signal

\be \label{1.1}
Z(t)=e^{i\vartheta(t)}\zf
\ee
that is generated by the Riemann zeta-function, where
\be \label{1.2}
\vartheta(t)=-\frac t2\ln\pi+\text{Im}\ln\Gamma\left(\frac 14+i\frac t2\right)=\frac t2\ln\frac{t}{2\pi}-\frac t2-\frac{\pi}{8}+\mcal{O}\left(\frac 1t\right) .
\ee
Let us remind that
\be \label{1.3}
\tilde{Z}^2(t)=\frac{{\rm d}\vp_1(t)}{{\rm d}t}, \ \vp_1(t)=\frac 12\vp(t)
\ee
where
\be \label{1.4}
\tilde{Z}^2(t)=\frac{Z^2(t)}{2\Phi^\prime_\vp[\vp(t)]}=\frac{Z^2(t)}{\left\{ 1+\mcal{O}\left(\frac{\ln\ln t}{\ln t}\right)\right\}\ln t}
\ee
(see \cite{12}, (5.1)-(5.3)) and $\vp_1(T),\ T\geq T_0[\vp_1]$ is the Jacob's ladder.

\subsection{}

It is known that the system of trigonometric functions
\be \label{1.5}
\left\{ 1, \cos\left(\frac{\pi}{l}t\right), \sin\left(\frac{\pi}{l}t\right),\dots, \cos\left(\frac{\pi}{l}nt\right),
\sin\left(\frac{\pi}{l}nt\right), \dots\right\}
\ee
is the orthogonal system on the segment $[0,2l]$. In this direction the following theorem holds true.

\begin{mydef11}
Let $\mcal{J}(2l)=\vp_1\{\mathring{\mcal{J}}(2l)\}$, where
\begin{eqnarray*}
& &
\mcal{J}(2l)=\mcal{J}(2l,K)=[2lK,2l(K+1)] , \\
& &
\mathring{\mcal{J}}(2l)=\mathring{\mcal{J}}(2l,K)=\left[\widering{2lK},\widering{2l(K+1)}\right], \ 2lK\geq T_0[\vp_1] , \\
& &
2l\in \left(\left. 0,\frac{T}{\ln T}\right.\right] ; \ K\in \mbb{N} .
\end{eqnarray*}
Then the system of functions
\be \label{1.6}
\left\{ 1, \cos\left(\frac{\pi}{l}\vp_1(t)\right), \sin\left(\frac{\pi}{l}\vp_1(t)\right),\dots ,
\cos\left(\frac{\pi}{l}n\vp_1(t)\right), \sin\left(\frac{\pi}{l}n\vp_1(t)\right),\dots \right\}
\ee
is the orthogonal system on $\mathring{\mcal{J}}(2l)$ with respect to the weight function $\tilde{Z}^2(t)$, i.e. the following new system of
integrals
\be \label{1.7}
\begin{split}
& \int_{\mathring{\mcal{J}}(2l)}\cos\left(\frac{\pi}{l}m\vp_1(t)\right)\cos\left(\frac{\pi}{l}n\vp_1(t)\right)\tilde{Z}^2(t){\rm d}t=
\left\{ \begin{array}{rcl} 0 & , & m\not= n , \\ l & , & m=n , \end{array} \right. \\
& \int_{\mathring{\mcal{J}}(2l)}\sin\left(\frac{\pi}{l}m\vp_1(t)\right)\sin\left(\frac{\pi}{l}n\vp_1(t)\right)\tilde{Z}^2(t){\rm d}t=
\left\{ \begin{array}{rcl} 0 & , & m\not= n , \\ l & , & m=n , \end{array} \right. \\
& \int_{\mathring{\mcal{J}}(2l)}\sin\left(\frac{\pi}{l}m\vp_1(t)\right)\cos\left(\frac{\pi}{l}n\vp_1(t)\right)\tilde{Z}^2(t){\rm d}t=0 , \\
& \int_{\mathring{\mcal{J}}(2l)}\cos\left(\frac{\pi}{l}n\vp_1(t)\right)\tilde{Z}^2(t){\rm d}t=0 , \\
& \int_{\mathring{\mcal{J}}(2l)}\sin\left(\frac{\pi}{l}n\vp_1(t)\right)\tilde{Z}^2(t){\rm d}t=0
\end{split}
\ee
for all $m,n\in\mbb{N}$ is obtained, where
\be \label{A}
t-\vp_1(t)\sim (1-c)\pi(t) , \tag{A}
\ee
\be \label{B}
2l(K+1)<\widering{2lK} , \tag{B}
\ee
\be \label{C}
\rho\{ \mcal{J}(2l);\mathring{\mcal{J}}(2l)\}\sim (1-c)\pi(t) \to\infty , \tag{C}
\ee
as $K\to\infty$, and $\rho$ denotes the distance of the corresponding segments, $c$ is the Euler constant and $\pi(t)$ is the prime-counting function.
\end{mydef11}

\begin{remark}
Theorem 1 gives the contact point between the functions $\zf,\ \pi(t),\ \vp_1(t)$ and the orthogonal system of trigonometric functions.
\end{remark}

\begin{remark}
It is clear that the formulae (\ref{1.7}) - for the modulated function $\tilde{Z}^2(t)$ - cannot be obtained in the known theories of Balasubramanian,
Heath-Brown and Ivic (comp. \cite{1}).
\end{remark}

This paper is a continuation of the series \cite{2}-\cite{15}.

\section{New method of the quantization of the Hardy-Littlewood integral (a special case)}

\subsection{}

We obtain from the first two formulae in (\ref{1.7})
\be \label{2.1}
\begin{split}
& \int_{\mathring{\mcal{J}}(2l)}\cos^2\left(\frac{\pi}{l}m\vp_1(t)\right)\tilde{Z}^2(t){\rm d}t=\frac 12|\mcal{J}(2l)| , \\
& \int_{\mathring{\mcal{J}}(2l)}\sin^2\left(\frac{\pi}{l}m\vp_1(t)\right)\tilde{Z}^2(t){\rm d}t=\frac 12|\mcal{J}(2l)|
\end{split}
\ee
for all $m\in\mbb{N}$. Next, from (\ref{2.1}) we obtain
\begin{cor}
\be \label{2.2}
\int_{\mathring{\mcal{J}}(2l)}\tilde{Z}^2(t){\rm d}t=|\mcal{J}(2l)|;\ |\mcal{J}(2l)|=2l .
\ee
\end{cor}

\subsection{}

Let us consider now the problem concerning the solid of revolution corresponding to the graph of the function (comp. \cite{5})
\bdis
\tilde{Z}(t),\ t\in [\widering{2lK},+\infty),\ 2lK>T_0[\vp_1] .
\edis

\begin{mydef3}
To divide this solid of revolution on parts of equal volumes.
\end{mydef3}

From (\ref{2.2}) we obtain the resolution of this problem.

\begin{cor}
Since
\be \label{a}
[\widering{2lK},+\infty)=\bigcup_{r=1}^\infty \mathring{\mcal{J}}(2l,r),\ \mathring{\mcal{J}}(2l,r)=
[\widering{2l(K+r-1)},\widering{2l(K+r)}] , \tag{a}
\ee
\be \label{b}
\pi\int_{\mathring{\mcal{J}}(2l,r)}\tilde{Z}^2(t){\rm d}t=2\pi l,\ r=1,2,3,\dots \ , \tag{b}
\ee
it follows that the sequence of points
\bdis
\{ \widering{2l(K+r-1)}\}_{r=2}^{+\infty}
\edis
is the resolution to the Problem for arbitrary fixed $2l\in (0,T/\ln T]$.
\end{cor}

\section{Generalization of the formula (\ref{2.2})}

\subsection{}

The following theorem holds true.

\begin{mydef12}
Let
\bdis
\mcal{J}(T,U)=[T,T+U],\ J(T,U)=\vp_1\{\mathring{\mcal{J}}(T,U)\};\ \mathring{\mcal{J}}(T,U)=[\mathring{T},\widering{T+U}] .
\edis
Then
\be \label{3.1}
\int_{\mathring{\mcal{J}}(T,U)}\tilde{Z}^2(t){\rm d}t=|\mcal{J}(T,U)|=U ,
\ee
for every $T\geq T_0[\vp_1],\ U\in (0,T/\ln T]$.
\end{mydef12}

\begin{remark}
From (\ref{3.1}) the general method for quantization of the Hardy-Littlewood integral follows (comp. Corollary 2:
$2lK\to \forall\ T\geq T_0[\vp_1],\ \mcal{J}(2l)\to\mcal{J}(T,U)$).
\end{remark}

Next, we obtain, using the mean-value theorem in (\ref{3.1})

\begin{cor}
\be \label{3.2}
\tilde{Z}^2(\xi)=\frac{|\mcal{J}(T,U)|}{|\mathring{\mcal{J}}(T,U)|},\ \xi\in \xi(\mathring{T},\widering{T+U}),\
\tilde{Z}(\xi)\not=0 ,
\ee
i.e.
\bdis
\tilde{Z}^2(\xi):1=|\mcal{J}(T,U)|:|\mathring{\mcal{J}}(T,U)| .
\edis
\end{cor}

\subsection{}

Let $\{[T',T'+1]\}$ stands for the continuum set of segments $[T',T'+1]\subset [T,T+T/\ln T]$. Since
\bdis
\frac{1}{|\mathring{\mcal{J}}(T',1)|}=\tilde{Z}^2(\xi),\ \xi=\xi(T')\in (\mathring{T}',\widering{T'+1})
\edis
then by the Riemann-Siegel formula
\bdis
Z(t)=2\sum_{n\leq\sqrt{\frac{t}{2\pi}}}\frac{1}{\sqrt{n}}\cos\{\vartheta(t)-t\ln n\}+\mcal{O}(t^{-1/4})
\edis
we obtain (see (\ref{1.4}))

\begin{cor}
\be \label{3.3}
\frac{1}{\sqrt{|\mathring{\mcal{J}}(T',1)|}}\sim\frac{2}{\sqrt{\ln \xi}}
\left|\sum_{n\leq\sqrt{\frac{\xi}{2\pi}}}\frac{1}{\sqrt{n}}\cos\{\vartheta(\xi)-\xi\ln n\}+\mcal{O}(\xi^{-1/4})\right|
\ee
where $\xi=\xi(T')$.
\end{cor}

\begin{remark}
The formula (\ref{3.3}) describes the complicated oscillations of the value $|\mathring{\mcal{J}}(T',1)|$ generated by the nonlinear transformation
$\mcal{J}(T',1)=\vp_1\{\mathring{\mcal{J}}(T',1)\}$.
\end{remark}

\section{Proof of Theorems 1 and 2}

\subsection{}

Let us remind that the following lemma is true (see \cite{12}, (5.1)-(5.3))

\begin{mydef5}
For every integrable function (in the Lebesgue sense) $f(x),\ x\in [\vp_1(T),\vp_1(T+U)]$ the following is true
\be \label{4.1}
\int_T^{T+U}f[\vp_1(t)]\tilde{Z}^2(t){\rm d}t=\int_{\vp_1(T)}^{\vp_1(T+U)}f(x){\rm d}x,\ U\in (0,T/\ln T] ,
\ee
where $t-\vp_1(t)\sim (1-c)\pi(t)$.
\end{mydef5}

\begin{remark}
The formula (\ref{4.1}) is true also in the case when the integral on the right-hand side of eq. (\ref{4.1}) is convergent but not absolutely (in the
Riemann sense).
\end{remark}

\subsection{}

If $\vp_1\{ [\mathring{T},\widering{T+U}]\}=[T,T+U]$ then we obtain from (\ref{4.1}) the following formula
\be \label{4.2}
\int_{\mathring{T}}^{\widering{T+U}}f[\vp_1(t)]\tilde{Z}^2(t){\rm d}t=\int_T^{T+U}f(x){\rm d}x,\ U\in (0,T/\ln T] .
\ee
Next, in the case $[T,T+U]=[2lK,2lK+2l]=\mcal{J}(2l)$, we have
\be \label{4.3}
\int_{\mcal{J}(2l)}F(t){\rm d}t=\int_0^{2l}F(t){\rm d}t
\ee
for every (integrable) $2l$-periodic function $F(t)$. Then from the known formulae

\bdis
\begin{split}
& \int_{\mcal{J}(2l)}\cos\left(\frac{\pi}{l}m\vp_1(t)\right)\cos\left(\frac{\pi}{l}n\vp_1(t)\right){\rm d}t=
\left\{ \begin{array}{rcl} 0 & , & m\not= n , \\ l & , & m=n , \end{array} \right. \\
& \int_{\mcal{J}(2l)}\sin\left(\frac{\pi}{l}m\vp_1(t)\right)\sin\left(\frac{\pi}{l}n\vp_1(t)\right){\rm d}t=
\left\{ \begin{array}{rcl} 0 & , & m\not= n , \\ l & , & m=n , \end{array} \right. \\
& \int_{\mcal{J}(2l)}\sin\left(\frac{\pi}{l}m\vp_1(t)\right)\cos\left(\frac{\pi}{l}n\vp_1(t)\right){\rm d}t=0 , \\
& \int_{\mcal{J}(2l)}\cos\left(\frac{\pi}{l}n\vp_1(t)\right){\rm d}t=0 , \quad
\int_{\mcal{J}(2l)}\sin\left(\frac{\pi}{l}n\vp_1(t)\right){\rm d}t=0, \quad m,n\in\mbb{N} ,
\end{split}
\edis
by the $\tilde{Z}^2$-transformation (see (\ref{4.2}), (\ref{4.3}); $[\mathring{T},\widering{T+U}]=\mathring{\mcal{J}}(2l)$) the
formulae (\ref{1.7}) follow. The properties (B), (C) in Theorem 1 are identical with \cite{13}, (A1), (B1).

\subsection{}

The formula (\ref{3.1}) follows from (\ref{4.2}) in the case $f(x)\equiv1$.

\section{Another type of the orthogonal systems}

It follows from (\ref{4.2}) that the continuum set $\mcal{S}(T,2l)$ of the systems
\bdis
\begin{split}
& \left\{ |\tilde{Z}(t)|, |\tilde{Z}(t)|\cos\left(\frac \pi l(\vp_1(t)-T)\right), |\tilde{Z}(t)|\sin\left(\frac \pi l(\vp_1(t)-T)\right),\dots , \right. \\
& \left.
|\tilde{Z}(t)|\cos\left(\frac \pi l n(\vp_1(t)-T)\right),|\tilde{Z}(t)|\sin\left(\frac \pi l n(\vp_1(t)-T)\right),\dots \right\} , \\
& t\in [ \mathring{T},\widering{T+2l}]
\end{split}
\edis
for all
\bdis
T\geq T_0[\vp_1],\ 2l\in (0,T/\ln T]
\edis
is the set of orthogonal systems on $[\mathring{T},\widering{T+2l}]$.

\begin{remark}
Let us call the elements of the system $\mcal{S}(T,2l)$ for fixed $2l\in (0,T/\ln T]$ and for all $T\geq T_0[\vp_1]$ as
\emph{the clones} of the known orthogonal trigonometric system
\bdis
\left\{ 1, \cos\left(\frac{\pi}{l}t\right), \sin\left(\frac{\pi}{l}t\right),\dots, \cos\left(\frac{\pi}{l}nt\right),
\sin\left(\frac{\pi}{l}nt\right), \dots\right\},\ t\in [0,2l] .
\edis
\end{remark}

\thanks{I would like to thank Michal Demetrian for helping me with the electronic version of this work.}

\end{document}